\pdfoutput=1
\documentclass[11pt,a4paper]{article}

\usepackage[utf8]{inputenc}
\usepackage[T1]{fontenc}
\usepackage{lmodern}
\usepackage[margin=1.05in]{geometry}
\usepackage{amsmath,amssymb,amsthm}
\usepackage{mathtools}
\usepackage{graphicx}
\usepackage{booktabs}
\usepackage{natbib}
\usepackage{microtype}
\usepackage{xcolor}
\usepackage[colorlinks=true,linkcolor=blue!50!black,citecolor=blue!50!black,urlcolor=blue!50!black]{hyperref}
\usepackage{enumitem}
\setlist{nosep}

\theoremstyle{definition}
\newtheorem{definition}{Definition}
\newtheorem{example}{Example}
\theoremstyle{plain}
\newtheorem{proposition}{Proposition}

\theoremstyle{remark}
\newtheorem{remark}{Remark}

\DeclareMathOperator{\sgn}{sgn}
\DeclareMathOperator{\Aut}{Aut}
\DeclareMathOperator{\id}{id}
\DeclareMathOperator{\Var}{Var}
\DeclareMathOperator{\corr}{corr}
\newcommand{\R}{\mathbb{R}}
\newcommand{\N}{\mathbb{N}}
\newcommand{\T}{\mathbb{T}}
\newcommand{\Z}{\mathbb{Z}}
\newcommand{\Prob}{\mathcal{P}}
\newcommand{\Multi}{\mathcal{M}}            
\newcommand{\Nbhd}{\mathcal{N}}             
\newcommand{\St}{\Sigma}                    
\newcommand{\Ty}{\mathcal{K}}               
\newcommand{\Env}{E}                        
\newcommand{\norm}[1]{\left\lVert #1 \right\rVert}
\newcommand{\ind}{\mathbf{1}}               
\newcommand{\mset}[1]{\{\!\{ #1 \}\!\}}
\newcommand{\Gen}{F}

\title{Generalized particle systems:\\ a locality-factorized metamodel for self-propelled particles and agent-based models}
\author{Igor Strozzi\\[2pt]
\small Applied Mathematics Department, Universidade Federal do Rio de Janeiro\\
\small \texttt{istrozzi@matematica.ufrj.br}}
\date{\today}

\begin{document}
\maketitle

\begin{abstract}
We propose a single formal object, the \emph{generalized particle system} (GPS), under which systems of self-propelled particles (SPPS) and agent-based models (ABMs) are instances of one construction. A GPS is a typed finite population, a per-particle state space of arbitrary internal structure, a global environment, a neighbourhood map that may depend on the current configuration, and type-indexed local update rules; these data assemble into a generator whose iteration is a monoid action in the sense of autonomous dynamical systems. Optional extensions---update schedules and stochastic local rules---are given at the same level of generality, so that sequential and noisy models (the original Primordial Particle System, the Vicsek model) are covered rather than approximated. The construction isolates the one ingredient that separates the two communities: in the classical ABM reading of \citet{Hinkelmann2011} the dependency graph is fixed, whereas in an SPPS it is a function of the state. We then introduce an \emph{exchangeability} axiom---rules depend on neighbours only through the multiset of their types and states---and prove that under it the generator is equivariant for type-preserving permutations and descends to a dynamical system on finite counting measures, which is the form in which kinetic and mean-field limits of swarming models are posed. We argue that the tuple should be read as a \emph{metamodel}: its structural components (state space, neighbourhood map, rule family) coordinatize a space of models, inside which the usual numerical parameters live. Six pre-registered numerical experiments, run on a single engine that is a literal implementation of the generator, test the claims: freezing the neighbourhood map destroys collective alignment ($\varphi = 0.91 \to 0.04$) while producing a different, non-aligned aggregation; the realized neighbourhood graph rewires at up to $64\%$ of edges per step and its rewiring rate is anti-correlated with order ($\corr = -0.99$); exchanging a metric for a topological neighbourhood measurably reshapes the phase diagram of a flocking rule; the same rule with cohesion and separation switched off exhibits a Vicsek-type order--disorder transition whose susceptibility peak grows with system size; and an internal activation variable does measurable but modest causal work. One pre-registered prediction failed and one subsumption claim could not be confirmed; both are reported.
\end{abstract}

\section{Introduction}

Two systems motivate this work. The \emph{Primordial Particle System} (PPS) of \citet{Schmickl2016, Schmickl2015} is a population of identical self-propelled particles in the plane whose heading update depends only on the number of neighbours to the left and to the right; at suitable densities it produces self-replicating, cell-like structures. \emph{Clusters} \citep{Ventrella2020}, and the closely related ``particle life'' systems \citep{CodeParade}, are populations of several particle \emph{types} interacting through a type-indexed force matrix; they produce orbiting, membrane-like and chasing structures. Both are deterministic or nearly so, both are local, and both exhibit non-trivial macroscopic organization despite the brevity of their defining rules. They belong to the family of \emph{self-propelled particle systems} (SPPS), whose canonical member is the Vicsek model \citep{Vicsek1995} and whose physics is reviewed in \citet{Vicsek2012, Marchetti2013}.

SPPS are rarely treated under the agent-based modelling (ABM) framework, although the particles are often called agents. Part of the reason is sociological---SPPS are studied mainly by physicists, ABMs by a broader community of complexity scientists---and part is technical: the ABM literature has produced formalizations in which the interaction structure is a fixed graph \citep{Hinkelmann2011, Mortveit2007}, whereas in an SPPS the interaction structure is decided geometrically at every instant. The present paper makes that difference precise, and shows that once it is isolated the two classes of models are instances of one construction.

\paragraph{Contributions.}
\begin{enumerate}[leftmargin=1.6em]
\item A definition (Definition~\ref{def:gps}) of a generalized particle system as a tuple $(I, \Ty, \kappa, \St, \Env, \Nbhd, \{f_k\}, g)$ whose generator $\Gen$ factors through purely local data. Fixed-graph ABMs, metric and topological SPPS, and typed particle-life systems are instantiations; update schedules and stochastic rules are added as extensions (Definitions~\ref{def:schedule}--\ref{def:stochastic}) without leaving the dynamical-systems setting.
\item An exchangeability axiom and three propositions: equivariance of $\Gen$ under type-preserving permutations (Proposition~\ref{prop:equivariance}); descent of the dynamics to the space of finite counting measures on $\Ty \times \St$ (Proposition~\ref{prop:descent}), which is the setting of kinetic limits of swarming models \citep{Cucker2007, HaTadmor2008, Carrillo2010, Degond2008}; and the identification of the fixed-neighbourhood case with (sequential) graph dynamical systems (Proposition~\ref{prop:gds}).
\item A metamodel reading of the definition, stated in a form that is operational: structural components of the tuple (discrete choices of $\St$, $\Nbhd$, $\{f_k\}$) are axes of a model space that are not redundant with the numerical parameters inside any one model.
\item Six pre-registered experiments (Section~\ref{sec:experiments}) on one engine that implements $\Gen$ literally, including two negative or unconfirmed outcomes.
\end{enumerate}

\paragraph{Related work.}
The formal lineage of our definition is that of graph dynamical systems and sequential dynamical systems \citep{Mortveit2007}, of which the ABM framework of \citet{Hinkelmann2011} is a polynomial-over-finite-field specialization; \citet{Smith2011} (``network automata'') and the adaptive-network literature \citep{Gross2008} study precisely the case in which the graph coevolves with the node states, which is the case we call state-dependent. On the modelling side, the ODD protocol \citep{Grimm2006, Grimm2010} and the DEVS formalism \citep{Zeigler2000} are specification standards for ABMs; our tuple is a mathematical rather than a descriptive specification and is compatible with either. Within SPPS, the distinction between metric and topological neighbourhoods is due to \citet{Ballerini2008}; the three flocking tendencies used in our running example are those of \citet{Reynolds1987}; weighted alignment with a kinetic limit is the Cucker--Smale model \citep{Cucker2007}. The programme of surveying \emph{spaces} of rules rather than parameters of one rule goes back to \citet{Wolfram1984} and \citet{Langton1990} for cellular automata, and has been pursued off-lattice in Lenia \citep{Chan2019}, its particle variant \citep{Mordvintsev2022}, and in interactive atlases of life-like rules \citep{Hovey2020}. Early work on multi-type self-propelled particles is \citet{Levine2001}; memory and internal energy depots as internal states appear in \citet{Nagai2015, Schweitzer1998}. The term ``interacting particle systems'' is taken in probability by the lattice models of \citet{Liggett1985}; our systems are off-lattice and (by default) deterministic, hence the different name.

\paragraph{Why a generalization.}
Beyond the usual benefits of unification---shared methods, and the isolation of the aspects that all such systems have in common---the specific payoff we pursue is the ability to survey \emph{model spaces}: to look, not for parameter regions in which a fixed model behaves interestingly, but for families of models that do. \citet{Hovey2020} maps the rule space of life-like cellular automata in exactly this sense; the systematic study of active-matter phases in \citet{Brilliantov2020} suggests the question is timely for particle systems too.

\section{Dynamical systems and the locality problem} \label{sec:ds}

\subsection{Autonomous dynamical systems}

\begin{definition}[Autonomous dynamical system] \label{def:ds}
Let $S$ be a set and $(T, \ast, e)$ a monoid. An \emph{autonomous dynamical system} on $S$ with time monoid $T$ is a left monoid action $\mathfrak{D} : T \times S \to S$, written $\mathfrak{D}_t(x) \coloneqq \mathfrak{D}(t, x)$, i.e.\ a map satisfying
\begin{align*}
	\mathfrak{D}_e &= \id_S, \\
	\mathfrak{D}_t \circ \mathfrak{D}_s &= \mathfrak{D}_{t \ast s} \qquad (s, t \in T).
\end{align*}
\end{definition}

When $T = (\N, +)$ the action is generated by the single map $\Gen \coloneqq \mathfrak{D}_1$ through $\mathfrak{D}_k = \Gen^{\circ k}$, and conversely any map $\Gen : S \to S$ generates one. All models below are of this discrete-time kind. We note for later that the same definition applies with $S$ replaced by a space of probability measures: a Markov kernel $P$ on $S$ generates the action $\mu \mapsto \mu P^k$ on $\Prob(S)$, so stochastic models are autonomous dynamical systems on $\Prob(S)$ \citep[cf.][for the random-dynamical-systems formulation]{Arnold1998}.

\subsection{A particle system, taken literally} \label{sec:pps}

Consider the PPS of \citet{Schmickl2016}. Particle $i$ has position $\mathbf{x}_i \in \R^2$ and heading $\theta_i \in \mathbb{S}^1$. With a fixed interaction radius $r$, let $N_i$ be the number of other particles within distance $r$ of $\mathbf{x}_i$ and let $L_i$, $R_i$ be the numbers of those in the left and right half-discs relative to the heading $\theta_i$. The update is
\begin{equation} \label{eq:pps}
	\begin{aligned}
		\theta_i' &= \theta_i + \alpha + \beta\, N_i\, \sgn(R_i - L_i), \\
		\mathbf{x}_i' &= \mathbf{x}_i + v\,(\cos\theta_i', \sin\theta_i'),
	\end{aligned}
\end{equation}
with $\alpha$ a fixed rotation, $\beta$ a turning gain and $v$ a fixed speed. (In the original model the particles are updated \emph{sequentially in a random order} within each time step; we return to this in Section~\ref{sec:schedule}. Equation~\eqref{eq:pps} is the synchronous version.)

As a dynamical system, \eqref{eq:pps} is a map $\Gen : (\R^2 \times \mathbb{S}^1)^n \to (\R^2 \times \mathbb{S}^1)^n$, and there is nothing problematic about that. What is worth making explicit is the \emph{shape} of its coordinate functions. Writing $s = (s_j)_{j \in I}$ for the configuration, the $i$-th coordinate of $\Gen(s)$ depends on $s_i$ and on $\{ s_j : \norm{\mathbf{x}_j - \mathbf{x}_i} \le r \}$---but \emph{which} $j$ enter is itself determined by $s$. Define the \emph{dependency graph at $s$} as the directed graph on $I$ with an edge $j \to i$ whenever $\Gen(s)_i$ depends on $s_j$. For \eqref{eq:pps} this is the geometric graph of radius $r$ on the current positions, and it changes at every step. The per-particle rule is simple; the coupling it induces on the configuration is a graph that is a function of the configuration.

\subsection{Agent-based models as dynamical systems} \label{sec:abm}

\citet{Hinkelmann2011} formalize an ABM as follows: given variables $\{x_i\}_{i \in I}$, $|I| = n$, taking values in a finite field (more generally a ring) $\mathbb{F}$, one fixes a directed \emph{dependency graph} $G$ on $I$ and, for each $i$, a local update function $f_i : \mathbb{F}^{\,\mathrm{in}(i)} \to \mathbb{F}$ depending only on the in-neighbours $\mathrm{in}(i)$ of $i$ in $G$. The coordinates assemble into $f = (f_1, \dots, f_n) : \mathbb{F}^n \to \mathbb{F}^n$, and a choice of \emph{update schedule} (parallel, sequential in a fixed order, etc.) turns $f$ into the generator of a discrete-time dynamical system. This is the specialization to polynomial functions over $\mathbb{F}$ of the graph dynamical systems and sequential dynamical systems of \citet{Mortveit2007}. Two features are to be retained: the dependency graph is \emph{fixed}, and the local functions are \emph{labelled}---$f_i$ may treat different in-neighbours differently (the Boolean rule $x_1' = x_2 \wedge \neg x_3$ is the paradigm case).

The contrast with Section~\ref{sec:pps} is now exact. A classical ABM has a fixed dependency graph and labelled local rules; an SPPS has a state-dependent dependency graph and, typically, rules that are symmetric in the neighbours. The next section defines an object that contains both.

\section{Generalized particle systems} \label{sec:gps}

\subsection{The definition}

Fix a finite index set $I = \{1, \dots, n\}$ of particles. For a set $X$, write $\Multi(X)$ for the set of finite multisets over $X$ (equivalently, finitely supported $\N$-valued measures on $X$), and $\Multi_m(X)$ for those of total mass $m$. For a finite set $J$ and a set $X$, write $X^J$ for the labelled families $(x_j)_{j \in J}$.

\begin{definition}[Generalized particle system] \label{def:gps}
A \emph{generalized particle system} (GPS) is a tuple
$$
	\mathcal{G} = \big(I,\ \Ty,\ \kappa,\ \St,\ \Env,\ \Nbhd,\ \{f_k\}_{k \in \Ty},\ g\big),
$$
where
\begin{itemize}[leftmargin=1.6em]
	\item $\Ty$ is a finite set of \emph{types} and $\kappa : I \to \Ty$ assigns a type to each particle;
	\item $\St$ is the \emph{per-particle state space}; a \emph{configuration} is an element $s = (s_i)_{i \in I} \in \St^I$;
	\item $\Env$ is the \emph{environment space}, carrying every global or non-particle quantity;
	\item $\Nbhd : I \times \St^I \times \Env \to \mathcal{P}(I)$ is the \emph{neighbourhood map}, with $\Nbhd(i, s, e) \subseteq I \setminus \{i\}$;
	\item for each $k \in \Ty$, the \emph{local rule} $f_k$ assigns to a particle's own state, the labelled family of its neighbours' types and states, and the environment, a new state:
	$$
		f_k : \St \times \Big( \bigsqcup_{J \subseteq I} (\Ty \times \St)^J \Big) \times \Env \longrightarrow \St;
	$$
	\item $g : \St^I \times \Env \to \Env$ is the \emph{environment update}.
\end{itemize}
These data assemble into the \emph{generator} $\Gen : \St^I \times \Env \to \St^I \times \Env$,
\begin{equation} \label{eq:generator}
	\Gen(s, e) = (s', e'), \qquad
	e' = g(s, e), \qquad
	s'_i = f_{\kappa(i)}\Big( s_i,\ \big( \kappa(j), s_j \big)_{j \in \Nbhd(i, s, e')},\ e' \Big),
\end{equation}
and the GPS is the autonomous dynamical system on $\St^I \times \Env$ with $T = (\N, +)$ generated by $\Gen$.
\end{definition}

The order of operations in \eqref{eq:generator}---environment first, then particles reading the updated environment---is a convention; the opposite order is obtained by composing $g$ with a one-step delay. When nothing global is needed one takes $\Env = \{\ast\}$ and omits $g$; we then write $\Nbhd(i, s)$ and $f_k(s_i, (\kappa(j), s_j)_{j \in \Nbhd(i,s)})$.

\begin{remark}[It is a dynamical system, and that is not the point]
Equation~\eqref{eq:generator} defines a single map, so $\mathfrak{D}_t = \Gen^{\circ t}$ satisfies the monoid-action axioms automatically. The content of Definition~\ref{def:gps} is not the dynamical-system structure, which is free, but the \emph{factorization} of the global map through local data $(\Nbhd, \{f_k\})$ and a separated global component $(\Env, g)$.
\end{remark}

\begin{remark}[The dependency graph]
The directed graph $G(s, e)$ on $I$ with edges $j \to i$ for $j \in \Nbhd(i, s, e)$ is the dependency graph of Section~\ref{sec:pps}, now an explicit component of the model. The classical ABM is the case $\Nbhd(i, s, e) = \mathrm{in}(i)$ independent of $(s, e)$; the SPPS case is a state-dependent $\Nbhd$, typically the \emph{metric} rule
\begin{equation} \label{eq:metric}
	\Nbhd_r(i, s) = \{\, j \ne i : d(x_i, x_j) \le r \,\},
\end{equation}
with $x_i$ the position component of $s_i$ and $d$ the metric of the ambient space (the flat torus $\T^2_L = (\R / L\Z)^2$ in all our examples), or the \emph{topological} rule $\Nbhd^{\mathrm{knn}}_k(i, s)$ = the $k$ nearest other particles \citep{Ballerini2008}. The two communities study the same object under different assumptions on this one ingredient.
\end{remark}

\begin{remark}[Environment and global state]
A rule that consults a running global statistic is \emph{not} an autonomous map on $\St^I$. The only way to keep such a rule inside Definition~\ref{def:gps} is to make the statistic a coordinate of $\Env$ and give it an update $g$. Example~\ref{ex:flock} needs this.
\end{remark}

\subsection{Extensions: schedules and noise} \label{sec:schedule}

Definition~\ref{def:gps} is synchronous and deterministic. Both restrictions are removed without leaving the setting of Definition~\ref{def:ds}.

For $i \in I$ let $\Gen_i : \St^I \times \Env \to \St^I \times \Env$ be the map that updates the $i$-th particle coordinate as in \eqref{eq:generator} (reading the environment $e' = g(s,e)$) and leaves all other coordinates fixed.

\begin{definition}[Schedule] \label{def:schedule}
A \emph{schedule} is a map $\sigma$ assigning to each $(s, e) \in \St^I \times \Env$ a finite word $\sigma(s, e) = (i_1, \dots, i_m)$ over $I$. The \emph{scheduled generator} is
$$
	\Gen_\sigma(s, e) = \big( \Gen_{i_m} \circ \cdots \circ \Gen_{i_1} \big)(s, e), \qquad (i_1, \dots, i_m) = \sigma(s, e),
$$
with the environment updated once, before the first particle. The synchronous generator \eqref{eq:generator} is the case in which all particles are updated from the same configuration $s$; a fixed ordering $\sigma \equiv (1, \dots, n)$ gives the sequential dynamical systems of \citet{Mortveit2007}.
\end{definition}

\begin{definition}[Stochastic GPS] \label{def:stochastic}
A \emph{stochastic} GPS replaces each local rule by a Markov kernel $f_k : \St \times \big( \bigsqcup_J (\Ty \times \St)^J \big) \times \Env \to \Prob(\St)$, the environment update by a kernel $g : \St^I \times \Env \to \Prob(\Env)$, and (optionally) the schedule by a kernel $\sigma : \St^I \times \Env \to \Prob(I^{<\omega})$. The synchronous stochastic generator is the kernel $\Gen$ on $\St^I \times \Env$ under which, given $(s, e)$, the environment is drawn from $g(s, e)$ and the new particle states are drawn \emph{conditionally independently} from $f_{\kappa(i)}(s_i, \cdot, e')$. The GPS is then the autonomous dynamical system $\mu \mapsto \mu \Gen^t$ on $\Prob(\St^I \times \Env)$; the deterministic case is recovered when all kernels are Dirac.
\end{definition}

The original PPS is a stochastic GPS with deterministic local rules and a uniformly random sequential schedule; the Vicsek model is a synchronous stochastic GPS with a deterministic alignment rule followed by angular noise (Section~\ref{sec:vicsek}).

\subsection{Exchangeability and what it buys} \label{sec:exchangeability}

The local rules of Definition~\ref{def:gps} are labelled, as in the ABM picture. Particle models almost never use the labels: a bird responds to the states of its neighbours, not to their names. We make this an axiom and derive its consequences.

\begin{definition}[Exchangeable GPS] \label{def:exchangeable}
A GPS is \emph{exchangeable} if
\begin{enumerate}[label=(E\arabic*),leftmargin=2.6em]
	\item \label{E1} each $f_k$ depends on the labelled family $(\kappa(j), s_j)_{j \in J}$ only through the multiset $\mset{ (\kappa(j), s_j) : j \in J } \in \Multi(\Ty \times \St)$; we then write $f_k : \St \times \Multi(\Ty \times \St) \times \Env \to \St$;
	\item \label{E2} the neighbourhood map and the environment update are equivariant under type-preserving relabellings: for every $\pi \in \Aut(\kappa) \coloneqq \{ \pi \in \mathrm{Sym}(I) : \kappa \circ \pi = \kappa \}$,
	$$
		\Nbhd(\pi i, \pi \cdot s, e) = \pi\, \Nbhd(i, s, e), \qquad g(\pi \cdot s, e) = g(s, e),
	$$
	where $(\pi \cdot s)_{\pi i} = s_i$.
\end{enumerate}
\end{definition}

Condition \ref{E2} holds for the metric rule \eqref{eq:metric}, for any $\Nbhd$ defined through a function of pairs of states, and for a constant $\Nbhd$ if and only if the fixed graph is $\Aut(\kappa)$-invariant. It holds for the $k$-nearest-neighbour rule up to tie-breaking (ties among equidistant particles must be broken by state, not by index, for \ref{E2} to hold exactly).

\begin{proposition}[Equivariance] \label{prop:equivariance}
If $\mathcal{G}$ is exchangeable then $\Gen(\pi \cdot s, e) = \pi \cdot \Gen(s, e)$ for every $\pi \in \Aut(\kappa)$, in the sense that $\Gen(\pi \cdot s, e) = (\pi \cdot s', e')$ whenever $\Gen(s, e) = (s', e')$.
\end{proposition}

\begin{proof}
By \ref{E2}, $g(\pi \cdot s, e) = g(s, e) = e'$. For the particle coordinate $\pi i$,
\begin{align*}
	\Gen(\pi \cdot s, e)_{\pi i}
	&= f_{\kappa(\pi i)}\Big( (\pi \cdot s)_{\pi i},\ \mset{ (\kappa(j'), (\pi\cdot s)_{j'}) : j' \in \Nbhd(\pi i, \pi \cdot s, e') },\ e' \Big) \\
	&= f_{\kappa(i)}\Big( s_i,\ \mset{ (\kappa(\pi j), s_j) : j \in \Nbhd(i, s, e') },\ e' \Big) \\
	&= f_{\kappa(i)}\Big( s_i,\ \mset{ (\kappa(j), s_j) : j \in \Nbhd(i, s, e') },\ e' \Big) = s'_i,
\end{align*}
using $\kappa \circ \pi = \kappa$, the substitution $j' = \pi j$ licensed by \ref{E2}, and \ref{E1} for the passage to multisets. Hence $\Gen(\pi \cdot s, e)_{\pi i} = s'_i = (\pi \cdot s')_{\pi i}$.
\end{proof}

Let $\rho : \St^I \to \Multi_n(\Ty \times \St)$, $\rho(s) = \mset{ (\kappa(i), s_i) : i \in I }$, be the \emph{typed empirical measure} of a configuration. Its fibres are exactly the $\Aut(\kappa)$-orbits.

\begin{proposition}[Descent to counting measures] \label{prop:descent}
If $\mathcal{G}$ is exchangeable there is a unique map $\bar{\Gen} : \Multi_n(\Ty \times \St) \times \Env \to \Multi_n(\Ty \times \St) \times \Env$ such that $(\rho \times \id_\Env) \circ \Gen = \bar{\Gen} \circ (\rho \times \id_\Env)$. Consequently the typed empirical measure $\mu_t = \rho(s_t)$ evolves autonomously: $(\mu_{t+1}, e_{t+1}) = \bar{\Gen}(\mu_t, e_t)$.
\end{proposition}

\begin{proof}
$\rho \times \id_\Env$ is surjective onto its image and its fibres are the orbits of $\Aut(\kappa)$ acting on the first factor. By Proposition~\ref{prop:equivariance}, $\Gen$ maps orbits to orbits, so $(\rho \times \id_\Env) \circ \Gen$ is constant on fibres of $\rho \times \id_\Env$ and factors uniquely through it.
\end{proof}

Proposition~\ref{prop:descent} is the formal reason the multiset formulation is the right one for \emph{kinetic} questions: a dynamics on empirical measures is precisely the starting point of the particle $\to$ kinetic $\to$ hydrodynamic hierarchy for swarming models \citep{Cucker2007, HaTadmor2008, Carrillo2010, Degond2008}. For the metric rule, $\bar{\Gen}$ has the explicit form
$$
	\bar\Gen(\mu, e)_{\text{particles}} = \sum_{(k, \sigma) \in \mathrm{supp}\,\mu} \mu(k, \sigma)\ \delta_{\big( k,\ f_k(\sigma,\ (\mu - \delta_{(k,\sigma)})|_{B_r(x(\sigma))},\ e') \big)},
$$
i.e.\ a particle in state $\sigma$ reads the restriction of the empirical measure (minus itself) to its own ball. A mean-field limit is then a statement about $\bar\Gen$ as $n \to \infty$ under a density scaling; we do not pursue it here, but note that Definition~\ref{def:exchangeable} is the hypothesis under which the question is well posed.

\begin{proposition}[Fixed neighbourhoods are graph dynamical systems] \label{prop:gds}
Let $\mathcal{G}$ be a (not necessarily exchangeable) GPS with $\Env = \{\ast\}$, $\Ty = I$, $\kappa = \id$, and $\Nbhd(i, s) = \mathrm{in}(i)$ for a fixed directed graph $G$ on $I$. Then the synchronous generator \eqref{eq:generator} is the parallel graph dynamical system $(G, \{f_i\})$ of \citet{Mortveit2007}, and the scheduled generator $\Gen_\sigma$ with $\sigma \equiv (i_1, \dots, i_n)$ a fixed permutation word is the sequential dynamical system with that update order. If moreover $\St = \mathbb{F}$ is a finite field and each $f_i$ is a polynomial function, this is the ABM of \citet{Hinkelmann2011}.
\end{proposition}

\begin{proof}
With $\Nbhd$ constant, \eqref{eq:generator} reads $s'_i = f_i(s_i, (j, s_j)_{j \in \mathrm{in}(i)})$, a function of the states of $i$ and its in-neighbours, which is the definition of a vertex function on $G$; the parallel and sequential compositions are those of the two cited frameworks by construction. The last claim is the definition in \citet{Hinkelmann2011}.
\end{proof}

Note that Proposition~\ref{prop:gds} needs the \emph{labelled} rules of Definition~\ref{def:gps}; under exchangeability one recovers only the symmetric graph dynamical systems. The hierarchy is therefore
$$
	\{\Nbhd \text{ constant}\} \ \subset\ \{\Nbhd \text{ a function of pairwise distances}\} \ \subset\ \{\Nbhd \text{ arbitrary}\},
$$
with classical ABMs, metric/topological SPPS, and adaptive networks \citep{Gross2008, Smith2011} as the three strata, and with the labelled/exchangeable distinction orthogonal to it.

\subsection{A metamodel, operationally} \label{sec:metamodel}

The tuple of Definition~\ref{def:gps} separates two kinds of choices. The numerical parameters of a model ($\alpha, \beta, r, v, \dots$) live \emph{inside} a fixed choice of $(\St, \Nbhd, \{f_k\})$; the latter data are the specification \emph{of} the model. Searching, within one model, for parameter regions of interesting behaviour is the usual practice, and it is the only practice available when the model is not written as a tuple.

We do not claim that ``the space of all $(\St, \Nbhd, \{f_k\})$'' is a mathematical object one can survey: it has no metric and no measure, and any concrete survey proceeds by fixing a parametrized family of rules, at which point the structural choices become additional---discrete, combinatorial---coordinates of a larger parameter space. The operational claim is weaker and testable: \emph{those structural coordinates are not redundant with the numerical ones}. Replacing a metric by a topological neighbourhood, freezing the neighbourhood map, switching a rule term off, or enlarging $\St$ by an internal variable, should move the system across the phase diagram in ways that no choice of the original numerical parameters reproduces. Section~\ref{sec:experiments} tests this for each of the four kinds of structural move.

\section{Instantiations} \label{sec:examples}

All examples live on the flat torus $\T^2_L$ with the torus metric in \eqref{eq:metric}; all have $T = (\N, +)$.

\begin{example}[The PPS, synchronous] \label{ex:pps}
Single type; $\St = \T^2_L \times \mathbb{S}^1$; $\Nbhd = \Nbhd_r$; $\Env = \{\ast\}$. From $s_i = (\mathbf{x}_i, \theta_i)$ and the neighbour multiset $\mset{(\mathbf{x}_j, \theta_j)}$ the rule computes $N_i$, the half-disc counts $L_i, R_i$ (which need the heading $\theta_i$ and the displacements $\mathbf{x}_j - \mathbf{x}_i$, and nothing else), and then \eqref{eq:pps}. The rule is exchangeable. The original model of \citet{Schmickl2016} is obtained from this by the uniformly random sequential schedule of Definition~\ref{def:schedule}; see Section~\ref{sec:pps-negative} for why the distinction is not cosmetic.
\end{example}

\begin{example}[An activation flock] \label{ex:flock}
Single type; $\St = \T^2_L \times \mathbb{S}^1 \times [0, 1]$ with state $s_i = (\mathbf{x}_i, \mathbf{v}_i, a_i)$, where $\mathbf{v}_i$ is a unit velocity and the scalar $a_i$ is an \emph{activation}, an internal variable beyond position and velocity. $\Nbhd = \Nbhd_r$; write $\Nbhd_i = \Nbhd_r(i, s)$, $n_i = |\Nbhd_i|$, $N_i = \max(n_i, 1)$. The velocity update superposes the three tendencies of \citet{Reynolds1987}---cohesion, separation, alignment---with cohesion \emph{weighted by the neighbours' activation}:
\begin{align*}
	\mathbf{c}_i &= \frac{f_c}{N_i} \sum_{j \in \Nbhd_i} (1 + a_j)\, (\mathbf{x}_j - \mathbf{x}_i), \qquad
	\mathbf{s}_i = \frac{f_s}{N_i} \sum_{\substack{j \in \Nbhd_i \\ \norm{\mathbf{x}_i - \mathbf{x}_j} < \rho}} (\mathbf{x}_i - \mathbf{x}_j), \qquad
	\mathbf{m}_i = \frac{f_m}{N_i} \sum_{j \in \Nbhd_i} \mathbf{v}_j, \\
	\mathbf{v}_i' &= \frac{\tilde{\mathbf{v}}_i}{\norm{\tilde{\mathbf{v}}_i}}, \quad
	\tilde{\mathbf{v}}_i = \tfrac12 \big( \mathbf{v}_i + \mathbf{c}_i + \mathbf{s}_i + \mathbf{m}_i \big), \qquad
	\mathbf{x}_i' = \mathbf{x}_i + V \mathbf{v}_i' \pmod L,
\end{align*}
where displacements are minimum-image displacements on the torus. The activation integrates the local density asymmetrically, with gates:
\begin{align*}
	a_i' = a_i + \eta_+\, \ind[\, n_i > \bar{n} \,]\, \ind[\, a_i < \gamma_+ \,] - \eta_-\, \ind[\, n_i \le \bar{n} \,]\, \ind[\, a_i > \gamma_- \,],
\end{align*}
which keeps $a_i \in [0, 1]$. A particle in a denser-than-typical neighbourhood becomes more active; an active particle exerts a stronger cohesive pull through $\mathbf{c}_i$. This is a positive feedback between local density and aggregation.

The threshold $\bar{n}$ is \emph{adaptive}: $\bar n = \tfrac12 (n_{\min} + n_{\max})$, where $n_{\min}, n_{\max}$ are the smallest and largest neighbourhood counts observed over the whole history of the run. As such the rule consults a global, history-dependent quantity, and is not an autonomous map on $\St^I$. It is brought inside Definition~\ref{def:gps} by $\Env = \N \times \N$, $e = (n_{\min}, n_{\max})$, and
$$
	g(s, e) = \Big( \min\big( e_1, \min_i n_i(s) \big),\ \max\big( e_2, \max_i n_i(s) \big) \Big),
$$
with the rule reading $\bar n$ from $e' = g(s, e)$. This is a minimal instance of the $(\Env, g)$ construction, and the one that forced it.

Reference parameters (used throughout Section~\ref{sec:experiments} unless stated): $r = 1.7$, $\rho = 0.8$, $V = 1.3$, $f_c = 0.25$, $f_s = 0.18$, $f_m = 0.084$, $\eta_+ = 0.05$, $\eta_- = 0.10$, $\gamma_+ = 0.94$, $\gamma_- = 0.11$; initial positions uniform on the torus, initial velocities uniform on the circle, $a_i(0) = 0$.
\end{example}

\begin{example}[Typed systems: \textit{Clusters} / particle life] \label{ex:plife}
Heterogeneity is the case $|\Ty| > 1$. Take $\St = \T^2_L \times \R^2$ (position, velocity), $\Nbhd = \Nbhd_r$, and a matrix $A \in \R^{\Ty \times \Ty}$, not necessarily symmetric. The rule of type $k$ reads the \emph{types} of its neighbours through the multiset in $\Multi(\Ty \times \St)$---this is why the second argument of $f_k$ carries types---and accumulates a force $\sum_j A_{k, \kappa(j)}\, \phi(\norm{\mathbf{x}_j - \mathbf{x}_i})\, \widehat{(\mathbf{x}_j - \mathbf{x}_i)}$ with a radial profile $\phi$ (repulsive below a core radius $r_{\min}$, then linear in $A$), followed by a damped velocity update. With $|\Ty| = 3$ this is the qualitative setting of \citet{CodeParade}; the discrete relatives of \citet{Levine2001} fit the same template. Types require no new structure: they are carried by $\kappa$, by $\{f_k\}$, and by the typed multiset.
\end{example}

\begin{example}[Vicsek as a structural move] \label{ex:vicsek}
In Example~\ref{ex:flock} set $f_c = f_s = 0$ and $f_m = 1$, and follow the velocity update by a rotation through an angle drawn uniformly from $[-\eta/2, \eta/2]$ (Definition~\ref{def:stochastic}). The result is a Vicsek-type alignment model \citep{Vicsek1995} in which a particle averages its own velocity with the mean velocity of its neighbours, i.e.\ with self-weight $1$ against neighbour-mean weight $1$---a variant of the standard rule, which weights self and each neighbour equally. Switching two rule terms off and adding noise is a change in $\{f_k\}$, not in a numerical parameter.
\end{example}

\section{Numerical experiments} \label{sec:experiments}

\subsection{Engine, protocol, and reproducibility}

All experiments were run on a single vectorized NumPy/SciPy engine whose central function is a literal implementation of \eqref{eq:generator}: a model is a \texttt{ModelSpec} holding $(\St, \Nbhd, \{f_k\}, \Env, g)$ as data, and one \texttt{step} function is the generator for every model. The engine reproduces, bit for bit, an earlier standalone implementation of Example~\ref{ex:flock}; this invariance is a permanent regression test, and every structural switch introduced below (neighbourhood type, freezing, rule ablation, noise) defaults to the original behaviour and is guarded by the same test. Neighbourhoods are computed with a periodic $k$-d tree. All runs are deterministic given $(\text{parameters}, \text{seed})$.

Each experiment was \emph{pre-registered} in the sense that its positive and negative criteria were written down before the run, and the governing filter was that a result counts as evidence for the framework only if it varies a \emph{structural} component of the tuple or measures something visible only under the dependency-graph view (a pure numerical sweep inside one model defends nothing about the generalization). Negative outcomes are reported. We state limitations up front: most runs use a single seed or two, system sizes are modest ($n \le 4096$), and the density regime for Experiments~\ref{sec:frozen}, \ref{sec:rewiring} and \ref{sec:activation} ($n = 1500$, $L = 50$) was selected by a density probe after a first pass at $n = 1000$, $L = 100$ proved too dilute to cluster at $r = 1.7$ (mean activation $\approx 0.1$); the model was not changed, only the regime.

Observables. For a configuration we use the \emph{polarization} $\varphi = \big\| \tfrac1n \sum_i \mathbf{v}_i \big\|$; the connected components of the undirected graph underlying $G(s)$ at the interaction radius, from which we take the \emph{giant fraction} (size of the largest component over $n$) and the \emph{mean cluster size}; and, for Experiment~\ref{sec:rewiring}, the \emph{rewiring rate}, the Jaccard distance between consecutive undirected edge sets, together with the mean degree and the transitivity of $G(s_t)$. Where a tail mean is reported it is over the last 40--120 steps as stated.

\subsection{Frozen versus state-dependent neighbourhood} \label{sec:frozen}

\emph{Claim tested:} Section~\ref{sec:pps}---the complexity of an SPPS is a state-dependent dependency graph, and the classical ABM is the fixed-graph special case. \emph{Design:} Example~\ref{ex:flock} with $n = 1500$, $L = 50$, 250 steps, seed 0, under two neighbourhood maps: the dynamic $\Nbhd_r$, and the \emph{frozen} map that fixes the edge set $(j \to i)$ at $t = 0$ and thereafter recomputes only the displacements from current positions (fixed partners, live geometry---the faithful fixed-graph reading of Proposition~\ref{prop:gds}). \emph{Pre-registered positive:} the frozen graph fails to develop the collective motion the dynamic graph develops. \emph{Negative:} it reproduces the macrostate.

\begin{table}[h]
\centering
\begin{tabular}{lcc}
\toprule
final macrostate & dynamic $\Nbhd$ (SPPS) & frozen $\Nbhd$ (ABM) \\
\midrule
polarization $\varphi$ & $0.905$ & $0.036$ \\
giant fraction & $0.151$ & $0.679$ \\
mean cluster size & $45.0$ & $74.1$ \\
\bottomrule
\end{tabular}
\caption{Experiment~\ref{sec:frozen}. RMS distance between the two trajectories: polarization $0.55$, giant fraction $0.32$.}
\label{tab:frozen}
\end{table}

\begin{figure}[h]
\centering
\includegraphics[width=\textwidth]{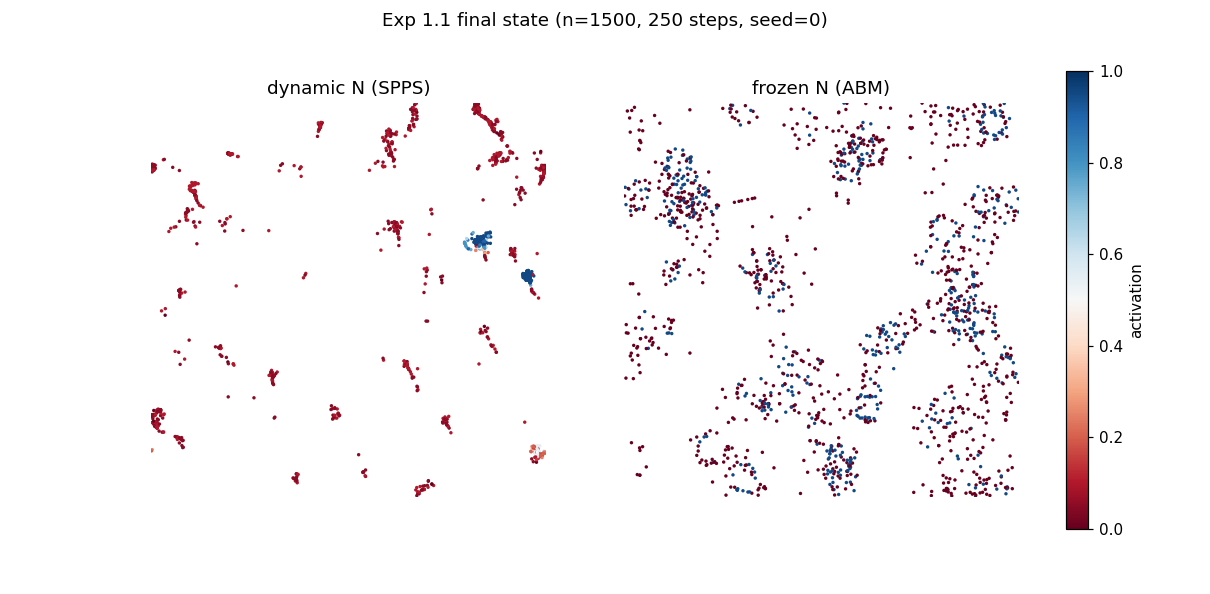}
\caption{Final configurations of Experiment~\ref{sec:frozen}, colour = activation. Left: dynamic neighbourhood, small aligned moving clusters, two of them activated. Right: frozen neighbourhood, large tangled non-aligned aggregates.}
\label{fig:frozen}
\end{figure}

\emph{Outcome: positive, with a nuance that should be kept.} The dynamic graph produces a coherent aligned state ($\varphi = 0.90$); the frozen graph does not align at all ($\varphi = 0.04$). But the frozen graph is not inert: with fixed partners and persistent cohesion it collapses into larger, tangled, non-aligned aggregates (higher giant fraction and mean cluster size). State-dependence of $\Nbhd$ is therefore what licenses \emph{alignment}; removing it does not remove aggregation but changes its kind. This sharpens the Section~\ref{sec:pps} claim: the macroscopic order of interest lives in the coupling, not in the per-particle rule, which is identical in both columns.

\subsection{The dependency graph as a temporal graph} \label{sec:rewiring}

\emph{Claim tested:} the state-dependent dependency graph is dynamically active, not a bookkeeping device. \emph{Design:} same run as Experiment~\ref{sec:frozen}, dynamic column; the realized graph $G(s_t)$ is recorded at every step and treated as a temporal graph. \emph{Pre-registered:} a non-trivial, time-varying rewiring rate and a measurable correlation between graph structure and $\varphi$; a near-zero or constant rewiring rate would mean the state-dependence is inert.

\begin{figure}[h]
\centering
\includegraphics[width=\textwidth]{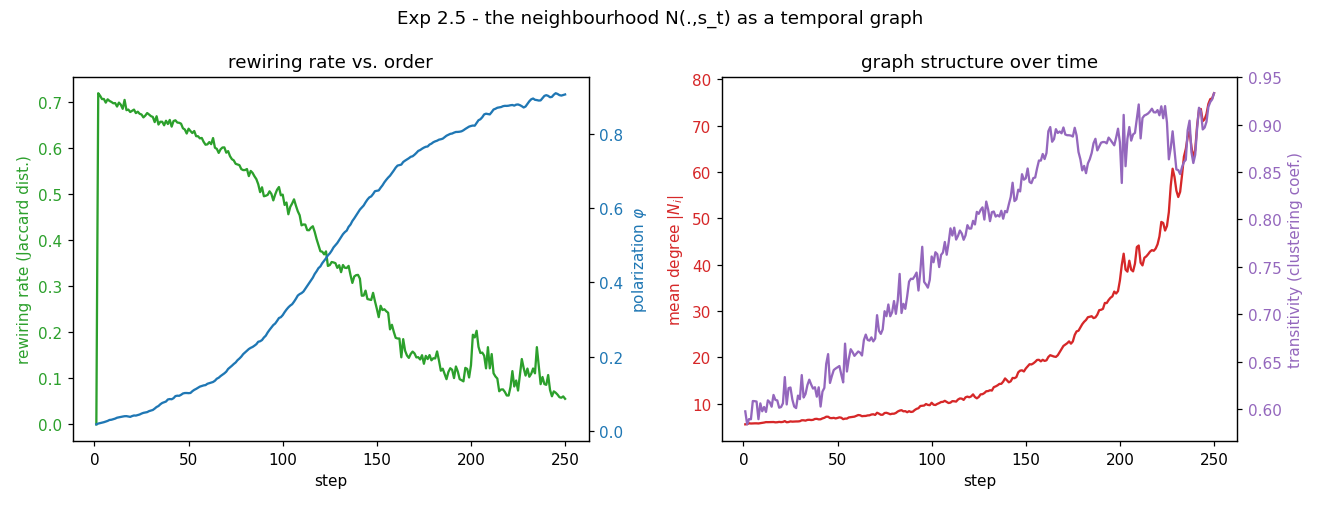}
\caption{Experiment~\ref{sec:rewiring}. Left: per-step rewiring rate (Jaccard distance of consecutive edge sets) against polarization. Right: mean degree and transitivity of $G(s_t)$.}
\label{fig:rewiring}
\end{figure}

\emph{Outcome: positive.} The rewiring rate starts near $0.7$ (two thirds of the edges change every step), averages $0.30$ over the run, and falls to $0.06$ as order sets in; the mean degree rises from $5.5$ to $77$ and the transitivity from $0.60$ to $0.93$. Over the post-warm-up window, $\corr(\text{rewiring rate}, \varphi) = -0.99$, $\corr(\text{transitivity}, \varphi) = +0.97$, $\corr(\text{mean degree}, \varphi) = +0.82$. The coupling structure progressively freezes as the flock orders. We stress that the near-perfect anti-correlation is a \emph{consistency} check rather than a discovery: aligned particles move together and therefore keep their neighbours. Its value is that it is a statement about the model that the per-particle treatment cannot formulate, and that it quantifies the rate at which an SPPS departs from its own fixed-graph approximation---here, by about a third of its edges per step in the disordered phase.

\subsection{A structural axis in the model-space atlas} \label{sec:atlas}

\emph{Claim tested:} Section~\ref{sec:metamodel}---a structural coordinate (the neighbourhood map) is not redundant with the numerical ones. \emph{Design:} Example~\ref{ex:flock} with $n = 800$, $L = 50$, 160 steps, tail 60, 2 seeds; a $6 \times 6$ grid in $(f_c, f_m) \in [0, 0.4] \times [0, 0.15]$ run twice, once with the metric $\Nbhd_r$ and once with the topological $\Nbhd^{\mathrm{knn}}_7$. \emph{Pre-registered positive:} structured regimes with boundaries, \emph{and} a measurable difference between the two panels. \emph{Negative:} featureless maps, or panels that coincide.

\begin{figure}[h]
\centering
\includegraphics[width=0.85\textwidth]{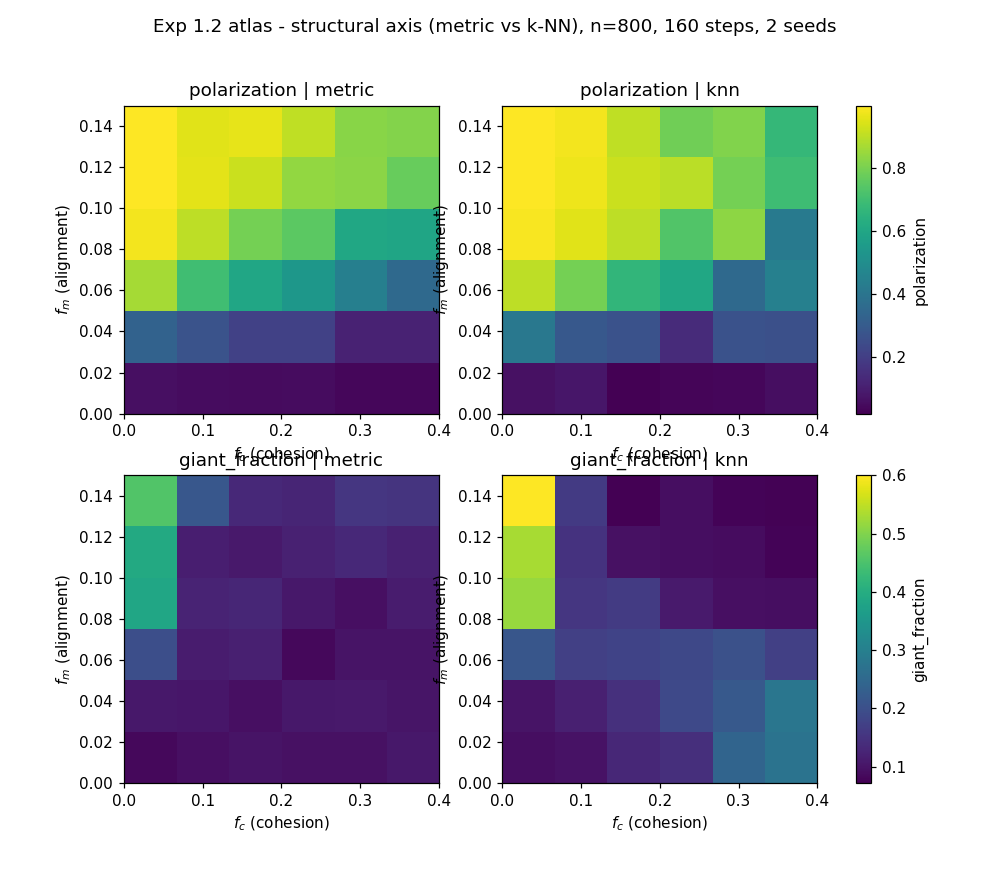}
\caption{Experiment~\ref{sec:atlas}. Polarization (top) and giant fraction (bottom) over $(f_c, f_m)$ for the metric (left) and $k$-nearest-neighbour (right) neighbourhood maps.}
\label{fig:atlas}
\end{figure}

\emph{Outcome: positive.} Polarization spans $[0.04, 0.99]$ (metric) and $[0.02, 1.00]$ ($k$-NN) across the plane, with a clear order--disorder boundary in $f_m$; the giant fraction spans $[0.08, 0.46]$ and $[0.07, 0.60]$ respectively. The RMS difference between the two panels is $0.082$ (polarization) and $0.080$ (giant fraction); the topological map pushes the giant fraction higher and shifts the aggregation region. The structural coordinate moves the system in a direction that is not available inside the numerical plane.

\subsection{A Vicsek-type transition from the flock rule} \label{sec:vicsek}

\emph{Claim tested:} subsumption of classical flocking by a move in $\{f_k\}$ (Example~\ref{ex:vicsek}). \emph{Design:} density $\rho = 0.66$, $r = 1.7$, $n \in \{256, 1024, 4096\}$ at fixed density ($L \approx 19.7, 39.4, 78.8$), 320 steps, tail 120, 2 seeds, 21 values of $\eta \in [0, 2\pi]$. \emph{Pre-registered positive:} $\varphi$ falls from $\approx 1$ to $\approx 0$ with $\eta$ and the transition sharpens with $n$. \emph{Negative:} no ordered phase, or no sharpening.

A first pass with a coarse $\eta$ grid used the maximal slope of $\varphi(\eta)$ as the sharpening diagnostic and returned the same, grid-limited slope for all three sizes; that diagnostic is retained in the configuration files but is not evidence either way. The reported run uses a finer grid and the susceptibility $\chi = n \Var(\varphi)$ over the tail.

\begin{table}[h]
\centering
\begin{tabular}{ccccc}
\toprule
$n$ & $\varphi(\eta = 0)$ & $\varphi(\eta = 2\pi)$ & $\chi$ at peak & $\eta$ at peak \\
\midrule
$256$  & $1.000$ & $0.054$ & $2.3$  & $2.51$ \\
$1024$ & $0.984$ & $0.028$ & $10.4$ & $2.51$ \\
$4096$ & $0.957$ & $0.014$ & $29.3$ & $1.88$ \\
\bottomrule
\end{tabular}
\caption{Experiment~\ref{sec:vicsek}.}
\label{tab:vicsek}
\end{table}

\begin{figure}[h]
\centering
\includegraphics[width=\textwidth]{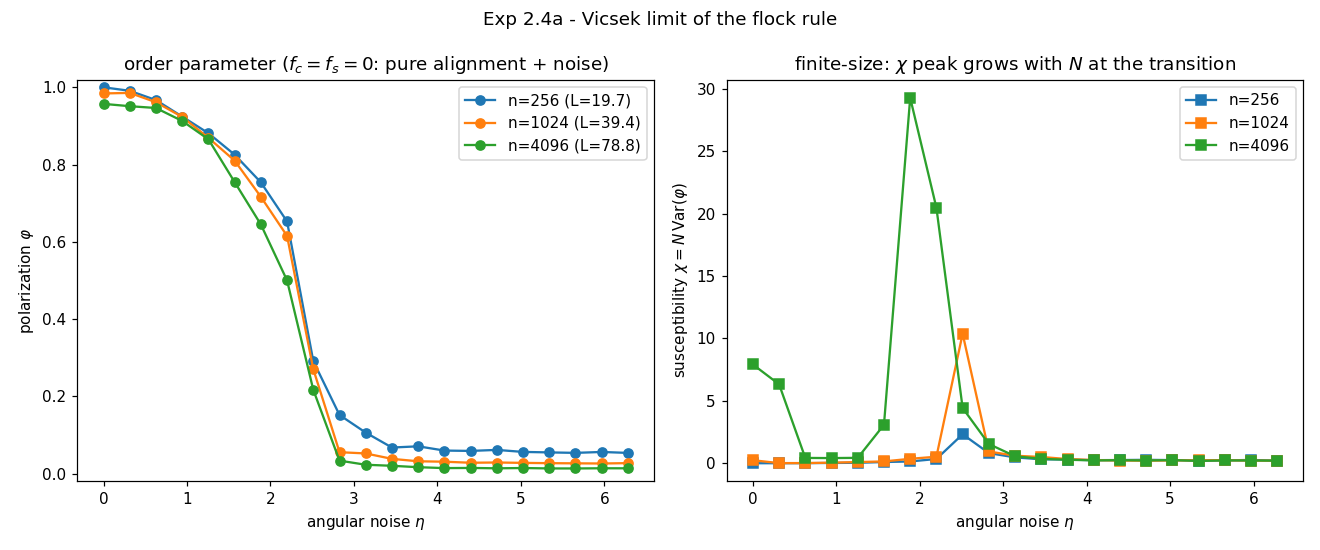}
\caption{Experiment~\ref{sec:vicsek}. Left: polarization against angular noise for three system sizes. Right: susceptibility $\chi = n\Var(\varphi)$.}
\label{fig:vicsek}
\end{figure}

\emph{Outcome: positive as pre-registered; the interpretation must be limited.} There is an ordered phase at low noise and a disordered one at high noise, with a crossover at $\eta_c \approx 1.9$--$2.5$, and the susceptibility peak grows with $n$ ($2.3 \to 10.4 \to 29.3$). This is the finite-size signature of a transition and confirms that the flock rule contains a Vicsek-type limit. We do \emph{not} claim to resolve the order of the transition: with $n \le 4096$ and two seeds the data cannot distinguish a continuous transition from the discontinuous one argued for the Vicsek model at large sizes \citep{Gregoire2004, Chate2008}, and the non-zero $\chi$ at $\eta = 0$ for $n = 4096$ shows that the ordered phase had not fully relaxed within 320 steps at that size. The peak heights are not linear in $n$ (ratios $4.5$ and $2.8$), so no scaling exponent is claimed.

\subsection{Does the internal state do causal work?} \label{sec:activation}

\emph{Claim tested:} the enlargement of $\St$ by an internal variable (Example~\ref{ex:flock}) is a structural move with consequences. \emph{Design (A, ablation):} $n = 1500$, $L = 50$, 250 steps, seed 0, tail 40; the weight $(1 + a_j)$ in $\mathbf{c}_i$ is replaced by $1$ while $a_i$ continues to evolve (so it is the \emph{coupling} that is ablated, not the variable). \emph{Design (B, robustness):} the adaptive threshold $\bar n$ is replaced by the per-step median, the per-step $0.75$-quantile, and an exponential moving average ($\alpha = 0.05$) of the per-step midpoint; the pre-registered prediction was that the running-extrema rule, being history-dependent and unprincipled, would be the most fragile component.

\begin{table}[h]
\centering
\begin{tabular}{lcc}
\toprule
(A) ablation & coupling on & coupling off \\
\midrule
giant fraction & $0.160$ & $0.128$ \\
mean cluster size & $39.2$ & $46.3$ \\
\bottomrule
\end{tabular}
\hspace{1.5em}
\begin{tabular}{lcccc}
\toprule
(B) threshold rule & running extrema & median & quantile$_{0.75}$ & EMA \\
\midrule
giant fraction & $0.160$ & $0.161$ & $0.124$ & $0.113$ \\
mean activation & $0.217$ & $0.496$ & $0.269$ & $0.433$ \\
\bottomrule
\end{tabular}
\caption{Experiment~\ref{sec:activation}.}
\label{tab:activation}
\end{table}

\emph{Outcome (A): positive but modest.} The feedback changes the macrostate (giant fraction $+26\%$ relative, mean cluster size $-15\%$): coupling concentrates a larger dominant cluster with smaller satellites. The effect is real and in the predicted direction, but not dramatic, consistent with activation firing only in the densest cores. \emph{Outcome (B): the pre-registered prediction failed.} The clustering macrostate is robust to the threshold rule (spread $0.05$ in giant fraction; running extrema and median are indistinguishable), although the mean activation itself is not. The history-dependent environment component is therefore \emph{not} a source of fragility for the clustering outcome in this regime, contrary to our expectation.

\subsection{One engine, three rule sets---and a subsumption not confirmed} \label{sec:pps-negative}

The engine was also used to run Examples~\ref{ex:flock}, \ref{ex:pps} and \ref{ex:plife} through the same \texttt{step} function, changing only the \texttt{ModelSpec}. The activation flock and the typed particle-life system behaved as expected (global alignment with a live internal variable in the first case, $\varphi = 0.67$; dense type-segregated clumps in the second, mean degree $105$). The \emph{synchronous} PPS (Example~\ref{ex:pps}; $n = 1200$, $L = 80$, $r = 5$, $v = 0.67$, $\alpha = \pi$, $\beta = 17^\circ$, 500 steps), however, did not form the cell-like structures reported by \citet{Schmickl2016}: the final configuration is visually indistinguishable from the uniform initial one, and the ``single percolating component'' returned by the clustering observable is the trivial percolation of a uniform configuration at that density ($\pi r^2 n / L^2 \approx 15$ expected neighbours), not a morphology. We report this as a subsumption claim \emph{not confirmed}. Two differences from the original are candidates: the density used here ($0.19$ particles per unit area) is above the range $0.04$--$0.08$ in which the cells are reported, and, more interestingly for this paper, the original updates particles sequentially in random order (Definition~\ref{def:schedule}) while our run is synchronous. The latter possibility is exactly the point of Section~\ref{sec:schedule}: the schedule is part of the model data, and ``the PPS'' as a synchronous deterministic GPS may simply be a different model from the PPS of \citet{Schmickl2016}. Resolving which difference is responsible is a well-posed experiment in the present framework (vary the schedule at fixed density, and the density at fixed schedule) and is left open.

\section{Discussion and conclusions}

We have given one definition under which the Primordial Particle System, Reynolds-type flocking, Vicsek-type alignment, and multi-type particle-life systems are instantiations, and which coincides, when the neighbourhood map is constant and the rules are labelled, with the graph-dynamical-system formulation of agent-based models. The device is modest: factor the global map through a neighbourhood map, type-indexed local rules, and a separated environment. It makes three things precise. The ``deceptive simplicity'' of particle rules is identified with a state-dependent dependency graph, and Experiments~\ref{sec:frozen}--\ref{sec:rewiring} show that this graph is where the collective order lives and that it is far from static. The exchangeability of particle rules, made an axiom, yields equivariance and descent to counting measures, which locates the framework at the starting point of the kinetic programme for swarming models. And the data specifying a model are separated from the parameters inside it, which, in the operational form of Section~\ref{sec:metamodel}, is a claim that Experiments~\ref{sec:atlas}--\ref{sec:activation} support for four kinds of structural move.

Limitations. The definition is discrete-time; the continuous-time analogue (a local vector field $v_k : \St \times \Multi(\Ty \times \St) \times \Env \to T\St$ with the same locality factorization) is natural but not developed, and most of active-matter theory lives there. The experiments are small (one or two seeds, $n \le 4096$) and were designed to be falsifiable rather than to be precise; none of the reported numbers should be read as an estimate with error bars. One subsumption claim (the sequential PPS) is unconfirmed, and one pre-registered prediction (fragility of the history-dependent threshold) failed.

Two directions seem natural. One is coarse-graining: aggregating a subset of particles into a single generalized particle---a quasi-particle with its own $\St$---so that ensembles can be described at a coarser level while remaining within Definition~\ref{def:gps}; Proposition~\ref{prop:descent} is a first step, since it already replaces configurations by typed counting measures. The other is the survey programme of Section~\ref{sec:metamodel}, for which the engine described in Section~\ref{sec:experiments} was built: instrumenting $\Gen$ so that one moves through structural choices of $(\St, \Nbhd, \{f_k\})$ and detects families of rules that produce emergent behaviour, in the manner of rule-space mapping for cellular automata.

\paragraph{Code and data.} The engine, the three model specifications, the six experiment scripts, their pre-registered criteria, and the configuration files from which every number in Section~\ref{sec:experiments} is read are available from the author.

\bibliographystyle{apalike}
\bibliography{gps_bibliography}

\end{document}